\documentclass{article}
\usepackage[utf8]{inputenc}
\usepackage{amsfonts}
\usepackage{amsthm}
\usepackage{amsmath}
\usepackage{xcolor}
\usepackage{hyperref}
\usepackage{mathtools}

\newcommand{\R}{\mathbb{R}}
\newcommand{\N}{\mathbb{N}}

\newcommand{\vcdim}{\mathrm{VC}\!-\!\mathrm{dim}}

\newtheorem{prop}{Proposition}

\newtheorem{thmx}{Theorem}

\newtheorem*{mainthm}{Main Theorem}

\theoremstyle{definition}

\theoremstyle{remark}
\newtheorem*{remark*}{Remark}

\usepackage{etoolbox}
\usepackage{hyperref}

\makeatletter
\newcounter{author}
\renewcommand*\author[1]{%
  \stepcounter{author}%
  \ifnum\c@author=1
    \gdef\@author{#1}%
  \else
    \xdef\@author{\unexpanded\expandafter{\@author\and#1}}%
  \fi
  \csgdef{author@\the\c@author}{#1}}
\newcommand*\email[1]{%
  \csgdef{email@\the\c@author}{#1}}
\newcommand*\address[1]{%
  \csgdef{address@\the\c@author}{#1}}
\AtEndDocument{%
  \xdef\author@count{\the\c@author}%
  \c@author=1
  \print@authors}
\newcommand*\print@authors{%
  \ifnum\c@author>\author@count
  \else
    \print@author{\the\c@author}%
    \advance\c@author by 1
    \expandafter\print@authors
  \fi}
\newcommand*\print@author[1]{%
  \par\medskip
  \begin{tabular}{@{}l@{}}%
    \textsc{Addresses of \csuse{author@#1}}\\
    \csuse{address@#1}\\
    \textit{E-mail address}:
    \href{mailto:\csuse{email@#1}}{\csuse{email@#1}}
  \end{tabular}}
\makeatother

\title{On the Vapnik-Chervonenkis dimension of products of intervals in $\R^d$}

\author{Alirio Gómez Gómez }
\address{Instituto de Ciência e Tecnologia da Universidade Federal de São Paulo
\\  Avenida Cesare Monsueto Giulio Lattes, 1201 -- Eugênio de Melo,\\
São José dos Campos, SP, Brazil\\
CEP 12247-014}
\email{gomezgomezalirio@gmail.com}

\author{Pedro L. Kaufmann  }
\address{Instituto de Ciência e Tecnologia da Universidade Federal de São Paulo
\\  Avenida Cesare Monsueto Giulio Lattes, 1201 -- Eugênio de Melo,\\
São José dos Campos, SP, Brazil\\
CEP 12247-014}
\email{plkaufmann@unifesp.br}
\date{}

\begin{document}
\maketitle

\begin{abstract}
    We study combinatorial complexity of certain classes of products of intervals in $\R^d$, from the point of view of Vapnik-Chervonenkis geometry. As a consequence of the obtained results, we conclude that the Vapnik-Chervonenkis dimension of the set of balls in $\ell_\infty^d$ -- which denotes $\R^d$ equipped with the sup norm -- equals $\lfloor (3d+1)/2\rfloor$. \let\thefootnote\relax\footnote{The second named author was supported by Grant 2016/25574-8, S\~ao Paulo Research Foundation (FAPESP)}
\end{abstract}

\section{Introduction}

A classifier $f$ on a measurable space $X$ is a binary function defined on $X$. A class $\mathcal{F}$ of classifiers on $X$ is said to  \textit{shatter} a sample $\sigma \subset X$ of size $n$ if $\mathcal{F}$ can perceive all possible binary labellings of the elements of $\sigma$, that is,  $$\# \{f|_{\sigma}: f\in \mathcal{F}\} =2^n.$$  In other  words, $\mathcal{F}$ shatters $\sigma$ when $\#\{C\cap \sigma : \chi_C \in \mathcal{F}\}=2^n$, where $\chi_C$ denotes the indicator function of $C$. Thus, we can alternatively focus on the sets $\{C: \chi_C \in \mathcal{F}\}$ to study properties of $\mathcal{F}$. More generally, families of classifiers on $X$ are in one-to-one correspondence with families of subsets of $X$. The later are called \emph{concept classes}, see for example \cite{pestov2019elementos}.
When  a learning  algorithm  has  to efficiently choose a  classifier within a family  $\mathcal{F}$ that minimizes the  learning error,
%(which means, intuitively, that such classifier that has been chosen  based on the data existing on a training set $\sigma$  will be efficient to predict the labelling of any point of $X$ outside of $\sigma$), 
it is often necessary some control on the quantity of different labellings that $\mathcal{F}$ can produce on finite samples of $X$.
%,  see \cite{pestov2019elementos}. 
One of the most important ways to measure the combinatorial complexity of families of classifiers, or equivalently, families of concept classes,  is to analyse their \emph{Vapnik-Chervonenkis dimension}, a concept introduced by Vapnik and Chervonenkis in \cite{vapnik-chervonenkis}.

The Vapnik-Chervonenkis dimension of a concept class $\mathcal{E}$ on $X$, which we will denote by $\vcdim(\mathcal{E})$, is defined by
$$
\vcdim(\mathcal{E})=\sup\{\#\sigma\mid\mbox{$\sigma\subset X$ is finite and $\mathcal{E}$ shatters $\sigma$}\}.
$$
We shall also write VC dimension, for short. To illustrate how the information about the VC dimension of a concept class $\mathcal{E}$ guarantees an efficient determination of an appropriate classifier, denote by $\mathcal{N}(\mathcal{E},n)$ the \emph{shattering coefficient} of $\mathcal{E}$ with respect to a sample size $n$. This is the number of labellings that $\mathcal{E}$ can produce on a sample of size $n$, that is, 
$$
\mathcal{N}(\mathcal{E},n)=\max_{\#\sigma=n} \#\{C\cap \sigma\mid C\in \mathcal{E} \}.
$$
%The number of  different outputs $Y_1, Y_2,..., Y_n$ that a function class $\mathcal{F}$ can achieve  on samples of  size $n$ is called the \textit{Shattering Coefficient of $\mathcal{F}$ with respect to the sample  size $n$  } and denoted by $\mathcal{N}(\mathcal{F},n)$. The  growth behavior of $\mathcal{N}(\mathcal{F},n)$ is characterized by $\vcdim(\mathcal{F})$. 
It was proved independently  by Sauer in 1972, Shelah in 1972, and Vapnik and  Chervonenkis in 1971 that if $\vcdim(\mathcal {E})=d<\infty$, then 
$\mathcal{N}(\mathcal{E},n)\leq (en/d)^d$,
where $e$ stands for the Euler constant. This implies in particular that the shattering coefficients grow polynomially with respect to the sample size. 
%Which implies that if $n\geq \vcdim(\mathcal{F})$ then shattering coefficients behaves like a polynomial functions of the  sample size $n$. 
Furthermore, empirical risk minimization is  consistent with respect to $\mathcal{E}$ if and only if $\vcdim(\mathcal{E})$ is finite, see \cite{lux2011statistical}.

The VC geometry of certain concept classes in $\R^n$ have received special attention. In his classic paper \cite{dudley1979balls}, Dudley showed that the VC dimension of Euclidean balls in $\R^n$ is $d+1$. Other natural classes that have been studied are the class of products of (possibly degenerate) intervals, 
\begin{eqnarray}\mathcal{R}:=\{\prod \limits_{i=1}^{n}[a_i, b_i]: -\infty \leq a_i< b_i\leq \infty\ , i=1,2,...,n\},
\label{eqWD}
\end{eqnarray}
and some subclasses of $\mathcal{R}$. In \cite{wenocur1981some}, Dudley and Wenocur proved that $\vcdim(\mathcal{R})=2n$. They showed additionally that, under the assumption that each $a_i$ in \eqref{eqWD} equals $-\infty$, the resulting subclass has VC dimension $n$. More recently, in \cite{gey2018vapnik}, Gey determined the VC dimension of the class of axis-parallel cuts in $\mathbb{R}^n$. More precisely, it  was proved that the VC dimension of the  concept class $\mathcal{A}_n=\{\{x=(x_1,x_2, \dots, x_n)\in \mathbb{R}^n: x_i\leq a \}: i=1,\dots, n, \hspace{0.2 cm} a \in \mathbb{R}\}$ is equal to $\max \{m: {m  \choose \lfloor m/2 \rfloor}\leq n \}$. 

\

In this  work we continue to investigate the VC complexity of some natural subclasses of $\mathcal{R}$. As a consequence of our study, we obtain in particular our main theorem, stated below. % we study the $\vcdim$ for  some class of subsets of $\mathbb{R}^n$ with sides parallel to coordinate axis. 
%The VC dimension  some classes of subsets of $\R^n$ has been widely  studied, for example Dudley and Wenocur, in  \cite{wenocur1981some}, proved that  the class $\mathcal{C}:=\{\prod \limits_{i=1}^{n}[a_i, b_i]: -\infty \leq a_i< b_i\leq \infty\ , i=1,2,...,n\}$ has VC dimension $2n+1$ and if the additional hypothesis $a_i=-\infty$ for all $i=1,2,...,n$ is assumed, then $\vcdim(\mathcal{C})=n+1$.  More recently, Gey in  \cite{gey2018vapnik} determined the VC dimension for the class of axis-parallel cuts in $\mathbb{R}^n$, more precisely it  was proved that the  concept class $\mathcal{A}_n=\{\{x=(x_1,x_2, \dots, x_n)\in \mathbb{R}^n: x_i\leq a \}: i=1,\dots, n, \hspace{0.2 cm} a \in \mathbb{R}\}$ has $\vcdim$ equals to $\max \{m: {m  \choose \lfloor m/2 \rfloor}\leq n \}$. In this  work  we study the $\vcdim$ for  some class of subsets of $\mathbb{R}^n$ with sides parallel to coordinate axis. 
Let us first establish some notation. $\ell_\infty^d$ denotes the vector space $\R^d$ equipped with the norm $\|\cdot\|_\infty:\mathbf{x}\to \max\{|x_1|,\dots,|x_d|\}$. Denote by $\mathcal{C}_d$ the set of all closed balls in $\ell_\infty^d$, which coincide with the set of closed cubes with sides parallel to the coordinate axes. Our main result reads as follows.
\begin{mainthm}
For each $d\geq 1$, $\vcdim(\mathcal{C}_d)=\lfloor (3d+1)/2\rfloor$. 
\end{mainthm}

This result was announced in \cite{despres2014vapnik}, but there was a flaw in the proof of the inequality $\vcdim(\mathcal{C}_d)\geq\lfloor (3d+1)/2\rfloor$. Specifically, in the proof of Lemma 2, it was claimed  that the set denoted by $S'$ is shattered by $\mathcal{C}_d$, which is false: the subset $\{\mathbf{x},\mathbf{y}\}$ cannot be carved out of $S'$ by an element of $\mathcal{C}_d$ (meaning that there is no $C\in\mathcal{C}_d$ such that $C\cap S'= \{\mathbf{x},\mathbf{y}\})$. This gap does not seem to be easily fixable.  

In the present work we prove the announced result approaching it from a different viewpoint. We study the Vapnik-Chervonenkis dimension of certain families of \emph{degenerate balls} in $\ell_\infty^d$. A closed (resp. open) degenerate ball in $\ell_\infty^d$ is any subset of the form $\prod_{i=1}^d I_i\subset \ell_\infty^d$, where each $I_i$ is a closed (resp. open) interval, unbounded to at least one side. As the notation suggests, a degenerate ball can be interpreted as a ball with infinite radius. For instance, if we choose a point $\mathbf{x}\in\ell_\infty^d$ and a direction $\mathbf{v}\in\ell_\infty^d\setminus\{\mathbf{0}\}$, and let $B_n$ be the closed ball centered at $n.\mathbf{v}$ and with $\mathbf{x}$ at its boundary, then the pointwise limit on $n\in\N$ of $B_n$ is a closed degenerate ball $D$ with $\mathbf{x}$ in its boundary. One can imagine that $D$ has its center at infinity, in the direction $\mathbf{v}$. Although the notation might be new, the idea was already used in Vapnik-Chervonenkis theory: in Dudley's proof in \cite{dudley1979balls} of the fact that the family of euclidean balls in $\R^d$ has VC dimension $d+1$, an underlying idea is that semi-spaces are pointwise limits of euclidean balls.  This motivates a more general and systematic study of combinatorial properties of degenerate balls in finite-dimensional normed spaces. 

Throughout this work, we shall use the following additional notation. For any closed subset $F$ of $\R^d$, we denote by $\mathcal{D}_d^F$ the set of all degenerate balls in $\ell_\infty^d$ containing $F$. 
%products of intervals $\prod_{i=1}^d I_i$ containing $F$ and such that each $I_i$ is a closed interval, unbounded to at least one side. 
In the case where $F=\{\mathbf{0}\}$, we use the lighter notation $\mathcal{D}_d^0$ instead of $\mathcal{D}_d^{\{\mathbf{0}\}}$. In the case where $F=\emptyset$, we simply write $\mathcal{D}_d$, instead of $\mathcal{D}_d^\emptyset$. Note that $\mathcal{D}_d$ is the set of all degenerate balls in $\ell_\infty^d$. Our first main result on Vapnik-Cervonenkis complexity of degenerate balls in $\ell_\infty^d$ is the following.

\begin{thmx}
For each nonempty bounded subset $F$ of $\R^d$, $\vcdim(\mathcal{D}^F_d)=\lfloor 3d/2\rfloor$. 
\end{thmx}

We then relate the Vapnik-Cervonenkis geometry of balls and of degenerate balls in $\ell_\infty^d$ spaces by proving the following. 

\begin{thmx}
For each $d\geq 2$, $\vcdim(\mathcal{C}_d)=\vcdim(\mathcal{D}^0_{d-1})+2$.
\end{thmx}

Note that Theorems A and B imply the Main Theorem, with the exception of  $\vcdim(\mathcal{C}_1)$ - but it is well known, and easily verified, that $\vcdim(\mathcal{C}_1)=2$. 

\

The next section  is dedicated to proving Theorems A and B. 

\section{Proof of main results}
 Let us start by establishing some notation. For each $S\subset \R^d$, $\mathbf{co}(S)$ denotes the convex hull of $S$, as usual. We define the \emph{rectangular hull of $S$} as being the smallest product of (possibly degenerate) intervals containing $S$, and denote it by $\Box$-hull($S$). 

To prove Theorems A and B, the following Propositions \ref{prop:DF} and \ref{prop:twitch} will be of use. 

%Consider the following families of subsets of $\R^d$: 
%\begin{align*}
 %   \mathcal{C}_d=&\{\mbox{cubes in }\R^d\},\\
  %  \mathcal{D}^F_d=&\{I=I_1\times\dots\times I_d:I_j\mbox{ of the form $(-\infty,a_j]$, $[a_j,\infty)$, or $\R$;}\mbox{ and }F\subset I\},\\
   % \mathcal{D}_d=&\mathcal{D}^{\{0\}}_d.
    %\mathcal{D}_d=\{&I_1\times\dots\times I_d:I_j\mbox{ closed intervals containing $0$ and with }0<|I_j|\leq\infty\}.
%\end{align*}

\begin{prop}\label{prop:DF} Let $F$ be a nonempty bounded closed subset of $\R^d$. 
Then, $\vcdim(\mathcal{D}^F_d)=\vcdim(\mathcal{D}^0_d)$.
\begin{proof}
Note that $\mathcal{D}_d^F=\mathcal{D}_d^R$, where $R=\Box$-hull($F$). We can assume then, without loss of generality, that $F$ is rectangular: $F=[a_1,b_1]\times\dots\times [a_d,b_d]$. For each $i\in\{1,\dots, d\}$, define
$$
p_i(x)=\begin{cases}
x-a_i,&\mbox{ if }x<a_i,\\
0,&\mbox{ if }a_i\leq x\leq b_i,\\
x-b_i,&\mbox{ if }b_i<x.
\end{cases}
$$
Define $P:\R^d\to \R^d$ by $P(p_1(x_1),\dots,p_d(x_d))$. It is readily verified that $P$ is surjective and satifies the following properties: 
\begin{enumerate}
    \item for each $D\in \mathcal{D}_d^F$, $P(D)\in \mathcal{D}^0_d$, and
    \item for each $D\in \mathcal{D}^0_d$, $P^{-1}(D)\in \mathcal{D}_d^F$. 
\end{enumerate}
Let us verify that $P$ also satisfies the following additional property: 
\begin{enumerate}
    \item[3.] $P$ is injective when restricted to any $\mathcal{D}_d^F$-shattered subset of $\R^d$. 
\end{enumerate}
Indeed, let $\mathbf{x}\in\R^d$, and let $D\in\mathcal{D}_d^F$ be such that $\mathbf{x}\in D$. Then 
$$
D\supset \mbox{$\Box$-hull}(F\cup\{\mathbf{x}\})
=\mathbf{co}([a_1,b_1]\cup\{x_1\})\times\dots\times \mathbf{co}([a_d,b_d]\cup\{x_d\}).
$$
Suppose now that $\mathbf{y}\in\R^d$ is such that $P(\mathbf{y})=P(\mathbf{x})$. Then for each $i\in\{1,\dots,d\}$, 
\begin{align*}
    p_i(\mathbf{y})=p_i(\mathbf{x})&
    \Rightarrow\begin{cases}y_i=x_i,&\mbox{ if }x_i<a_i\mbox{ or }x_i> b_i\\
    y_i\in [a_i,b_i],&\mbox{ if }x_i\in [a_i,b_i]
    \end{cases}\\
    &\Rightarrow y_i\in \mathbf{co}([a_i,b_i]\cup\{x_i\}). 
\end{align*}
It follows that $\mathbf{y}\in \Box$-hull$(F\cup\{\mathbf{x}\})\subset D$. This proves that $\mathbf{x}$ and $\mathbf{y}$ cannot be separated by an element of $\mathcal{D}_d^F$, from which Property 3 follows. \medskip

Let $S\subset\R^d$ be a $\mathcal{D}_d^F$-shattered set. By Property 3, $\#P(S)=\#S$. Let us verify that $P(S)$ is $\mathcal{D}^0_d$-shattered. Indeed, let $P(S')$ be an arbitrary subset of $P(S)$. Since $S$ is $\mathcal{D}_d^F$-shattered, there exists $D\in\mathcal{D}_d^F$ with $D\cap S=S'$. By Property 1, $P(S)\in\mathcal{D}_d^0$, and on the other hand, $P(D)\cap P(S)=P(S')$. This shows that $P(S)$ is $\mathcal{D}^0_d$-shattered. We conclude that $\dim(\mathcal{D}_d^0)\geq \mathcal{D}_d^F$. 

Suppose now that $S\subset \R^d$ is $\mathcal{D}_d^0$-shattered. For each $\mathbf{s}\in S$, choose any $\mathbf{t}_\mathbf{s}\in P^{-1}(\mathbf{s})$, and define $T=\{\mathbf{t}_\mathbf{s}:\mathbf{s}\in S\}$. We claim that $T$ is  $\mathcal{D}_d^F$-shattered. Indeed, let $T'=\{\mathbf{t}_\mathbf{s}:\mathbf{s}\in S'\}$ be any subset of $T$. Since $S$ is $\mathcal{D}_d^0$-shattered, there exists $D\in \mathcal{D}_d^0$ such that $S\cap D=S'$. By Property 2, $P^{-1}(S)\in \mathcal{D}_d^F$, and it clearly satisfies $P^{-1}(S)\cap T=T'$. This shows that $\dim(\mathcal{D}_d^F)\geq \mathcal{D}_d^0$, and we are done. \end{proof}
\end{prop}

\begin{prop}\label{prop:twitch} Let $\mathcal{B}$ be $\mathcal{D}_d^F$ or $\mathcal{C}_d$. Suppose that there exists a $\mathcal{B}$-shattered set $A\subset \R^d$ with $\# A=n$. Then, there exists a $\mathcal{B}$-shattered set $A'\subset \R$ with $\# A'=n$ and such that, for each coordinate projection $\pi_j:\R^d\to\R$, $\# \pi_j(A')=n$. 

\begin{proof} We shall prove the statement for $\mathcal{B}=\mathcal{C}_d$. The proof for $\mathcal{D}_d^F$ can be easily adapted. Suppose that $A\subset \R^d$ is $\mathcal{C}_d$-shattered and that $\# A=n$. Let $\mathbf{x}$ be some point in $A$. Consider $C_1,\dots,C_k\in \mathcal{C}_d$ that shatter $A$. For each $j=1,\dots,k$ such that $\mathbf{x}\in C_j$, note that $d(A\setminus\{\mathbf{x}\},C_j)>0$, since $C_j$ is closed. This implies that we can substitute $C_j$ by another cube $C'_j$ with same center but slightly bigger, so that $C'_j\cap A=C_j\cap A$, but now we guarantee that $\mathbf{x}$ is an interior point of $C'_j$. Choose some open set $V_j$ with $\mathbf{x}\in V_j \subset C'_j$. For the $j=1,\dots,k$, such that $\mathbf{x}\not\in C_j$, put $C'_j=C_j$ and choose an open neighborhood $V_j$ of $\mathbf{x}$ such that $V_j\cap C'_j=\emptyset$. 

Let $V=\cap_{j=1}^kV_j$. Note that $C'_1,\dots,C'_k$ shatter $A$, but also shatter $(A\setminus\{\mathbf{x}\})\cup\{\mathbf{x}'\}$ for any choice of $\mathbf{x}'\in V$. Since the set
$$
\cup_{j=1}^d \pi_j^{-1}(\pi_j(A\setminus\{\mathbf{x}\}))
$$
is nowhere dense in $\R^d$, we can choose some $\mathbf{x}'\in V$ such that, for each $j=1,\dots,d$, $\pi_j(\mathbf{x})\not\in\pi_j(A\setminus\{\mathbf{x}\})$. 

Repeating this process recursively to each point of $A$ at a time, we obtain $A'$, which clearly satisfies the desired properties.
\end{proof}
\end{prop}

\begin{remark*} The same works for for the set of non-degenerate rectangles in $\R^d$, or the set of closed balls with respect to any norm in $\R^d$. These cases will not be used in what follows, though. 
\end{remark*}

We are in position to prove the main results. 

\begin{proof}[Proof of Theorem A]
Let us start by showing that  $\vcdim(\mathcal{D}_d^0)\geq \lfloor 3d/2\rfloor$. Suppose that $d$ is even. In this case, we can write $3d/2$ instead of $\lfloor 3d/2\rfloor$.  The proof in this case will follow by a two-step induction on $d$. Note that $\mathcal{D}_2^0\geq 3$, since for instance it is readily verified that  $\{(-1,1),(1,-1),(2,1)\}$ is $\mathcal{D}_2^0$-shattered. Suppose now that $\vcdim(\mathcal{D}_{d}^0)\geq 3d/2$. This means that there exists a $\mathcal{D}_{d}^0$-shattered set  $S\subset \R^d$ with $\#S=3d/2$. We shall show that there is a $\mathcal{D}_{d+2}^0$-shattered subset of $\R^{d+2}$ with $3(d+2)/2=\#S+3$ elements, implying that $\vcdim(\mathcal{D}_{d+2}^0)\geq 3(d+2)/2$. Indeed, let $X$ %$\{\mathbf{x},\mathbf{y},\mathbf{z}\}$ 
be a $\mathcal{D}_2^0$-shattered subset of $\R^2$ with $\#X=3$. We claim that the set 
$$
T=\{(\mathbf{s},\mathbf{0}):\mathbf{s}\in S\}\cup\{(\mathbf{0},\mathbf{x}):\mathbf{x}\in X\}
%,(\mathbf{0},\mathbf{y}),(\mathbf{0},\mathbf{z})\}
\subset \R^d\times \R^2=\R^{d+2}
$$
is $\mathcal{D}_{d+2}^0$-shattered. In effect, let $T'$ be any subset of $T$. Then there are $S'\subset S$ and $X'\subset X$ with $T'=\{(\mathbf{s},\mathbf{0}):\mathbf{s}\in S'\}\cup\{(\mathbf{0},\mathbf{x}):\mathbf{x}\in X'\}$. Since $S$ is $\mathcal{D}_d^0$-shattered, there exists $I=\prod_{i=1}^d I_i\in \mathcal{D}_d^0$ with $S\cap I=S'$. Since $X$ is $\mathcal{D}_2^0$-shattered, there exists $J=\prod_{i=d+1}^{d+2} I_i\in \mathcal{D}_2^0$ with $X\cap J=X'$. Then, $I\times J=\prod_{i=1}^{d+2} I_j\in \mathcal{D}_{d+2}^0$ satisfies $T\cap I\times J= T'$. Since $\#T=3(d+2)/2$, we are done for $d$ even. 
% Fix some $0<\epsilon<1$. Denote by $\{e_1,\dots,e_d\}$ the canonical basis of $\R^d$, and for each $j\in\{1,\dots,d/2\}$, define $a_j=\sum_{i=1}^d\lambda_{ij}e_i$, where
% $$
% \lambda_{ij}=\begin{cases}
% -1,&\mbox{ if }i\in \{2j-1,2j\}, \mbox{ and}\\
% \epsilon,&\mbox{ otherwise}.
% \end{cases}
% $$
% We claim that $S=\{e_1,\dots,e_d,a_1,\dots,a_{d/2}\}$ is $\mathcal{D}_d$-shattered, which will conclude the proof. In effect, let $S'$ be an arbitrary subset of $S$. For each $i\in \{1,\dots,d\}$, let $j$ be the unique integer such that $i\in\{2j-1,2j\}$, and define
% $$
% I_i=\begin{cases}
% [-1,\infty),&\mbox{ if }e_i\in S'\mbox{ and } a_j\in S',\\
% [0,\infty),&\mbox{ if }e_i\in S'\mbox{ and } a_j\not\in S',\\
% (-\infty,\epsilon],&\mbox{ if }e_i\not\in S'\mbox{ and } a_j\in S',\\
% (-\infty,0],&\mbox{ if }e_i\not\in S'\mbox{ and } a_j\not\in S'. 
% \end{cases}
% $$
% Direct verification shows that $I=\prod_{i=1}^d I_i\in\mathcal{D}_d$ satisfies $I\cap S=S'$. It follows that $\vcdim(\mathcal{D}_d)\geq \#S=3d/2=\lfloor 3d/2\rfloor$, for $d$ even. \smallskip

Now suppose that $d$ is odd. It is a general observation that $\vcdim(\mathcal{D}_d^0)>\vcdim(\mathcal{D}_{d-1}^0)$, since whenever $S\subset \R^{d-1}$ is shattered by $\mathcal{D}_{d-1}^0$, then the set $\{(\mathbf{x},0)\mid \mathbf{x}\in S\}\cup \{(\mathbf{0},1)\}\subset \R^d$ is clearly shattered by $\mathcal{D}_{d}^0$. It follows that, also for odd $d$, 
$$\vcdim(\mathcal{D}_d^0)\geq\vcdim(\mathcal{D}_{d-1}^0)+1=\frac{3(d-1)}{2}+1=\frac{3d}{2}-\frac12=\lfloor \frac{3d}{2}\rfloor.$$

\medskip

It remains to  show that $\vcdim(\mathcal{D}_d^0)\leq \lfloor 3d/2\rfloor$, for each $d\geq 1$. This part of the proof departs from the main idea used in \cite[Theorem 1]{despres2014vapnik}. Suppose that $S$ is a $\mathcal{D}_d^0$-shattered subset of $\R^d$, with $\#S=d+n$. 
%Denote by $[a_1,b_1]\times\dots\times [a_d,b_d]$ the rectangular envelope of $S$. 
For each $i\in\{1,\dots,d\}$, choose $\mathbf{l}_i$ and $\mathbf{u}_i$ with minimal, and respectively maximal, $i$th coordinate amongst $S$. % in $S$ such that $\pi_i(\mathbf{l}_i)=a_i$ and $\pi_i(\mathbf{u}_i)=b_i$. 
Note that each $\mathbf{s}\in S$ should appear on the list 
\begin{eqnarray}\label{eq:list}
\mathbf{l}_1,\mathbf{u}_1,\mathbf{l}_2,\mathbf{u}_2,\dots,\mathbf{l}_d,\mathbf{u}_d. 
\end{eqnarray}
Indeed, it would be otherwise impossible to carve out $S\setminus\{\mathbf{s}\}$ from $S$ with an element of $\mathcal{D}_d^0$, since $\mathbf{s}$ would be in the rectangular envelope of $S\setminus\{\mathbf{s}\}$. 

Let $k$ be the number of elements of $S$ which appear on the list \eqref{eq:list} exactly once. Assume that $k\geq d+1$. By the pigeonhole principle, there exists an $i$ such that both $\mathbf{l}_i$ and $\mathbf{u}_i$ appear on the list exactly once. Without loss of generality, assume that $i=d$. Let $[a_1,b_1]\times\dots\times [a_d,b_d]$ the the rectangular envelope of $S\setminus \{\mathbf{l}_d,\mathbf{u}_d\}$. 
%Let $\pi:\R^d\to \R^{d-1}$ be the projection onto the first $d-1$ coordinates. 
Since $\mathbf{l}_d$ and $\mathbf{u}_d$ appear on \eqref{eq:list} only once, we have that   
\begin{eqnarray}\label{eq:liui}
a_i\leq \pi_i(\mathbf{l}_d)\leq b_i,\mbox{ and } a_i\leq \pi_i(\mathbf{u}_d)\leq b_i, \mbox{ for each }i\in\{1,\dots,d-1\}.
\end{eqnarray} Suppose that a product of intervals $I=\prod_{i=1}^dI_i$ carves out $S\setminus \{\mathbf{l}_d,\mathbf{u}_d\}$ from $S$. $I$ must contain $[a_1,b_1]\times\dots\times [a_d,b_d]$. It follows from \eqref{eq:liui} and the fact that $\mathbf{l}_d$ has minimal $d$th coordinate among $S$, that $I_d$ is bounded from below by $\pi_d(\mathbf{l}_d)$. Analogously, $I_d$ is bounded from above by $\pi_d(\mathbf{u}_d)$. $I_d$ is therefore a bounded interval, and it follows that $I$ is not an element from $\mathcal{D}^0_d$. This shows in particular that  $S\setminus \{\mathbf{l}_d,\mathbf{u}_d\}$ cannot be carved out from $S$ by an element from $\mathcal{D}^0_d$, a contradiction. Then, $k{\color{blue}\leq} d$. 

Once the $k$ points of $S$ that appear only once in \eqref{eq:list}, $2d-k$ slots remain to be filled with the $d+n-k$ points of $S$, which appear on \eqref{eq:list} at least twice. It follows that $2(d+n-k)\leq 2d-k$, thus $2n{\color{red}\leq} k$. We conclude that 
$$
\#S=d+n{\color{red}\leq} d+\frac{k}{2}{\color{blue}\leq} \frac{3d}{2}.
$$

%In particular, $\pi(\mathbf{l}_d)$ is contained in the rectangular envelope of $\pi(S\setminus\{\mathbf{l}_d,\mathbf{u}_d\})$ in $\R^{d-1}$, which we denote by $R$. Analogously, $\pi(\mathbf{d}_d)$ is contained in $R$. We claim that $S\setminus \{\mathbf{l}_d,\mathbf{u}_d\}$ cannot be carved out from $S$ by an element of $\mathcal{D}_d$. In effect, suppose that $D=\prod_{i=1}^d\in\mathcal{D}_d$ contains $S\setminus \{\mathbf{l}_d,\mathbf{u}_d\}$. Since Since $\mathbf{l}_d$ appears on \eqref{eq:list} only once, for each coordinate $i\in\{1,\dots,d-1\}$
%$ =\min \pi_i(S)$ and $\mathbf{u}_i=\max \pi_i(S)$. 
\end{proof}

\begin{proof}[Proof of Theorem B] Let $A\subset \R^d$ be a $\mathcal{C}_d$-shattered set with $\#A=\vcdim(\mathcal{C}_d)$. By Proposition \ref{prop:twitch}, we can assume that $\# \pi_j(A)=\# A$, for each $j=1,\dots,d$. Let $I_1\times\dots\times I_d$ be the rectangular envelope of $A$. Without loss of generality, assume that $|I_d|=\max\{|I_1|,\dots,|I_d|\}$, where $|I_j |$ denotes the length of the  interval $I_j$. 

Write $I_d=[a_d,b_d]$, and let $\mathbf{x},\mathbf{y}\in A$ be such that $\pi_d(\mathbf{x})=a_d$ and $\pi_d(\mathbf{y})=b_d$. Let $A'$ be some subset of $A\setminus \{\mathbf{x},\mathbf{y}\}$. Since $A$ is $\mathcal{C}_d$-shattered, there must be some $C_{A'}\in\mathcal{C}_d$ with $C_{A'}\cap A=A'\cup\{\mathbf{x},\mathbf{y}\}$. Consider $\pi:\R^d\to \R^{d-1}$ the projection onto the first $d-1$ coordinates, and note that $\pi(C_{A'})$ is a $(d-1)$-dimensional cube such that $\pi(C_{A'})\cap \pi(A) = \pi(A'\cup\{\mathbf{x},\mathbf{y}\})$. Since $C_{A'}$ contains $\mathbf{x}$ and $\mathbf{y}$, the side of the cube $\pi(C_{A'})$ is greater than the diameter of $\pi(A)$ in $\ell_\infty^{d-1}$. It follows that there exists $D_{A'}\in\mathcal{D}^{\{\mathbf{x},\mathbf{y}\}}_{d-1}$ such that $D_{A'}\cap \pi(A) = \pi(A'\cup\{\mathbf{x},\mathbf{y}\})$. It follows that $\{D_{A'}:A'\subset A\setminus\{\mathbf{x},\mathbf{y}\}\}\subset \mathcal{D}^{\{\mathbf{x},\mathbf{y}\}}_{d-1}$ shatters the set $\pi(A\setminus\{\mathbf{x},\mathbf{y}\})$ in $\R^{d-1}$. Since $\#\pi(A\setminus\{\mathbf{x},\mathbf{y}\})=\#A-2=\vcdim(\mathcal{C}_d)-2$, it follows that $\vcdim(\mathcal{D}^{\{\mathbf{x},\mathbf{y}\}}_{d-1})\geq \vcdim(\mathcal{C}_d)-2$. From Proposition \ref{prop:DF} it follows that $\vcdim(\mathcal{D}^0_{d-1})\geq \vcdim(\mathcal{C}_d)-2$.\medskip 

Now let $A\subset \R^{d-1}$ be finite, $\mathcal{D}^0_{d-1}$-shattered set, and let $L>0$ be such that $A$ is contained in the closed $\ell_\infty^d$ ball centered at $\mathbf{0}$ and with radius $L/2$.
%and denote by $L$ the diameter of $A$ with respect to $\|\cdot\|_\infty$. 
Consider in $\R^d$ the points $\mathbf{x}=(0,...,0,L)$ and $\mathbf{y}=(0,...,0,-L)$, and let
$A'=\{(\mathbf{a},0):\mathbf{a}\in A\}\cup \{\mathbf{x},\mathbf{y}\}$. To show that $\vcdim(\mathcal{C}_{d})\geq \vcdim(\mathcal{D}^0_{d-1})+2$, it suffices to verify that $A'$ is $\mathcal{C}_{d}$-shattered.  In effect, let $B \subset A'$. Since $\pi(B)\setminus \{\mathbf{0}\}\subset A$, $\pi(B)\setminus \{\mathbf{0}\}$ is carved out by some $D \in \mathcal{D}^0_{d-1}$. Note that, for each $M\geq L$, $\pi(B)\setminus \{\mathbf{0}\}$ can also be carved out by some $C \in \mathcal{C}_{d-1}$ which contains $\mathbf{0}$ and has $\|\cdot\|_\infty$-diameter $M$. Let us define an element $C'\in \mathcal{C}_d$ depending on $B\cap\{\mathbf{x},\mathbf{y}\}$, as follows:
\begin{enumerate}
    \item if $B\cap\{\mathbf{x},\mathbf{y}\}=\{\mathbf{x},\mathbf{y}\}$, consider some $C \in\mathcal{C}_{d-1}$ which contains $\mathbf{0}$ and has $\|\cdot\|_\infty$-diameter $2L$. Define $C'=C\times [-L,L]$; 
    \item if $B\cap\{\mathbf{x},\mathbf{y}\}=\{\mathbf{x}\}$, consider some $C \in\mathcal{C}_{d-1}$ which contains $\mathbf{0}$ and has $\|\cdot\|_\infty$-diameter $L$. Define $C'=C\times [0,L]$; 
    \item analogously, if $B\cap\{\mathbf{x},\mathbf{y}\}=\{\mathbf{y}\}$, consider some $C \in\mathcal{C}_{d-1}$ which contains $\mathbf{0}$ and has $\|\cdot\|_\infty$-diameter $L$. Define $C'=C\times [-L,0]$; 
    \item if $B\cap\{\mathbf{x},\mathbf{y}\}=\emptyset$, consider some $C \in\mathcal{C}_{d-1}$ which contains $\mathbf{0}$ and has $\|\cdot\|_\infty$-diameter $L$. Define $C'=C\times [-L/2,L/2]$; 
\end{enumerate}
In each case, $C'$ carves $B$ out of $A'$, which concludes the proof. \end{proof}

% {\color{red}Since $B$ is finite and $\pi(B\setminus \{\mathbf{x},\mathbf{y}\})\subset A$ is carved out by some $D \in \mathcal{D}^0_{d-1}$,
% %and $\mathbf{0} \in D$, 
% we can find $D' \in \mathcal{C}_{d-1}$ such that $D'\subset D$ and $D'$ carves out $\pi (B \setminus \{\mathbf{x,y}\})$. Now, we  choose $I_d$ as follows: if $\{\mathbf{x},\mathbf{y}\}\setminus B = \{\mathbf{x},\mathbf{y}\}$ we put $I_d=[-L,L]$; if $\{\mathbf{x},\mathbf{y}\}\setminus B=\{\mathbf{x}\}$, we put $I_d=[-L-1,0]$; if $\{\mathbf{x},\mathbf{y}\}\setminus B=\{\mathbf{y}\}$, we put $I_d=[0,L+1]$, and if $\{\mathbf{x},\mathbf{y}\}\setminus B=\emptyset$, we put $I_d=[-L-1,L+1]$. 
% %equal to  $[-d_\infty$-diam($B$),$d_\infty$-diam($B$)] $(-\infty, 0]$, $[0, \infty)$, or $\mathbb{R}$, depending on the set $\{x,y\} \setminus B$ being $\{x,y\}$, $\{x\}$, $\{y\}$, or $\emptyset$ respectively. 
% Then $D'\times I_d$ carves out $B$. It follows that $\vcdim(\mathcal{C}_{d})\geq \vcdim(\mathcal{D}^0_{d-1})+2$. (o conjunto $D'\times I_d$ não necessariamente é um cubo)} 

\section{Final remarks}

The natural follow up to this work would be to determine the exact Vapnik-Chervonenkis dimension of the set $\mathcal{D}_d$ of all degenerate balls in $\ell_\infty^d$. Note that one can easily obtain the comparison
\begin{eqnarray*}
\lfloor \frac{3d}{2}\rfloor=\vcdim(\mathcal{D}_d^0)\leq \vcdim(\mathcal{D}_d) \leq \vcdim(\mathcal{C}_d)=\lfloor \frac{3(d+1)}{2}\rfloor.
\label{floor}
\end{eqnarray*}
Indeed, the first inequality is clear since $\mathcal{D}_d^0\subset \mathcal{D}_d$. The second inequality follows from the following fact: if we can carve out a subset $S'$ from a set $S$ with a degenerate ball, then we can carve out $S'$ from $S$ with an appropriate ball with large enough radius. Computing $\lfloor 3d/2\rfloor$ and $\lfloor (3d+1)/2\rfloor$ in the cases where $d$ is odd and even separately gives us the following result, which is a direct consequence of Theorem A and the main theorem. 

\begin{prop}
Let $d\geq 1$. If $d$ is odd, $\vcdim(\mathcal{D}_d)=(3d+1)/2.$
If $d$ is even, 
$$
\frac{3d}{2}\leq \vcdim(\mathcal{D}_d)\leq \frac{3d}{2}+1. 
$$
\end{prop}

\bibliographystyle{plain}
\bibliography{references}

\begin{thebibliography}{1}

\bibitem{despres2014vapnik}
C.~J.~J. Després.
\newblock The {V}apnik-{C}hervonenkis dimension of cubes in $\mathbb{R}^d$.
\newblock {\em arXiv preprint arXiv:1412.6612}, 2017.

\bibitem{dudley1979balls}
R.~M. Dudley.
\newblock Balls in $\mathbb{R}^k$ do not cut all subsets of $k+ 2$ points.
\newblock {\em Advances in Mathematics}, 31(3):306--308, 1979.

\bibitem{gey2018vapnik}
S.~Gey.
\newblock Vapnik-{C}hervonenkis dimension of axis-parallel cuts.
\newblock {\em Communications in Statistics - Theory and Methods},
  47(9):2291--2296, 2018.

\bibitem{pestov2019elementos}
V.~G. Pestov.
\newblock {\em Elementos da {T}eoria de {A}prendizagem de {M}\'aquina
  {S}upervisionada}.
\newblock To appear in the 32$^{nd}$ ed. of the Brazilian Quolloquia of
  Mathematics Series, IMPA. arXiv preprint arXiv:1910.06820 (in Portuguese),
  2019.

\bibitem{lux2011statistical}
U.~V.~Luxburg and B.~Sch{\"o}lkopf.
\newblock Statistical learning theory: Models, concepts, and results.
\newblock In {\em Handbook of the History of Logic}, volume~10, pages 651--706.
  North-Holland, 2011.

\bibitem{vapnik-chervonenkis}
V.N. Vapnik and A.~Ya. Chervonenkis.
\newblock On the uniform convergence of relative frequencies of events io their
  probabilities.
\newblock {\em Theor. Probability Appl.}, 16:264--280, 1971.

\bibitem{wenocur1981some}
R.~S. Wenocur and R.~M. Dudley.
\newblock Some special {V}apnik-{C}hervonenkis classes.
\newblock {\em Discrete Mathematics}, 33(3):313--318, 1981.

\end{thebibliography}

%%%%%%%%%%%%

\end{document}